\documentclass[reqno,10pt]{amsart} 

\usepackage[colorlinks = true, linkcolor=blue,
            urlcolor=red,
            citecolor=olive]{hyperref}
\usepackage[utf8]{inputenc}
\usepackage[T1]{fontenc}
\usepackage{amsmath}

\usepackage{mathtools}      
\usepackage{mathabx}        

\usepackage{tikz}
\usepackage{tikz-3dplot}

\usepackage{amssymb}
\usepackage{amsfonts}
\usepackage{graphicx}

\usepackage{multicol}
\usepackage{subfiles} 

\usepackage[colorinlistoftodos]{todonotes}

\hypersetup{bookmarksdepth  =  3} 

\newtheorem{thm}{Theorem}[section]
\newtheorem{lem}[thm]{Lemma}
\newtheorem{prop}[thm]{Proposition}

\newtheorem{defn}[thm]{Definition}

\newtheorem{rem}[thm]{Remark}

\newtheorem{que}[thm]{Question}

\newcommand{\R}{\mathbb{R}}
\newcommand{\C}{\mathbb{C}}

\newcommand{\Id}{\text{Id}}

\DeclareMathOperator{\vis}{vis}
\newcommand{\partialvis}{\partial_{\vis}}

\DeclareMathOperator{\End}{End}

\DeclareMathOperator{\tr}{tr}

\DeclareMathOperator{\SL}{SL}
\DeclareMathOperator{\GL}{GL}
\DeclareMathOperator{\Hom}{Hom}

\DeclareMathOperator{\SO}{SO}

\makeatletter
\DeclareRobustCommand*{\bfseries}{%
  \not@math@alphabet\bfseries\mathbf
  \fontseries\bfdefault\selectfont
  \boldmath
}
\makeatother

\newcommand{\ie}{i.e. }

\newcommand{\calE}{\mathcal{E}}
\newcommand{\calK}{\mathcal{K}}
\newcommand{\calO}{\mathcal{O}}
\newcommand{\dbar}{\overline\partial}

\title{Anosov deformations of Barbot representations}

\copyrightinfo{2024}{Samuel Bronstein, Colin Davalo}

\author{Colin Davalo}
\address{Dipartimento di Matematica "Giuseppe Peano"
, Università di Torino, Via Carlo Alberto, 10 - 10123 Torino, Italy}
\email{colinmarcthierry@unito.it}

\author{Samuel Bronstein}
\address{Max Plank Institute for Mathematics in the Sciences, Inselstraße 22 - 04103 Leipzig, Germany}
\email{samuel.bronstein@mis.mpg.de}

\date{\today}
\begin{document}

\begin{abstract}
We construct for each conformal structure on a closed orientable surface of genus at least $2$ a proper slice in the character variety of representations of the associated surface group into $\SL(3,\R)$ that belongs to the Barbot component and show that the corresponding representations are Borel Anosov. We describe a fibered geometric structure in the space of full flags associated with these representations.

\end{abstract}
\maketitle

\tableofcontents
\section{Introduction.}

Let $S_g$ be a closed orientable surface and $\Gamma_g$ be its fundamental group. Among representations $\rho:\Gamma_g\to \SL(3,\R)$, an exceptional connected component was studied by Hitchin \cite{Hit87}. This connected component contains the composition of any Fuchsian representation with the irreducible representation $\iota_{\text{irr}}:\SL(2,\R)\to\SL(3,\R)$. Every representation in this connected component is discrete and faithful, and even \emph{Anosov} \cite{CG93}, which is a stronger dynamical property introduced by Labourie \cite{Labourie} that amounts to having large eigenvalues gaps for the images of the representations. Moreover these representations are the holonomies of convex projective structures on the surface.

\medskip

The other two connected components of the character variety contain representations with compact image \cite{Hit87}, which therefore cannot be discrete and faithful. However one of these two components also contains some Fuchsian representations that are Anosov. Indeed one can compose any Fuchsian representation with the reducible representation $\iota_{\text{red}}:\SL(2,\R)\to\SL(3,\R)$. Barbot studied these representations and certain deformations inside $\GL(2,\R)\subset \SL(3,\R)$ \cite{Bar01},\cite{Bar10}. For such deformations, which are all discrete and faithful, he identifies exactly at which point they stop being Anosov. He also considers some other deformations in the subgroup of affine transformations of $\R^2$ inside $\SL(3,\R)$, these representations are not anymore semi-simple representations and hence do not correspond to new points in the character variety.

\medskip

In this paper we consider proper slices of the character variety containing only Anosov representations, which are necessarily discrete and faithful. These slices intersect the deformations studied by Barbot only at the Fuchsian locus.

\medskip

More precisely, given the choice of a Riemann surface structure $X$ on $S_g$ and the choice of a square root of the canonical bundle $\mathcal{K}$ on $X$, we construct a properly embedded disk in the character variety of (real) dimension $4g-4$ parametrized by elements of $H^0(\mathcal{K}^{\frac{3}{2}})$. Note that the expected dimension of the character variety is $16g-16$, and the space of possible Riemann surface structures, the Teichmüller space of $S_g$, has dimension $6g-6$. 

\medskip

Given a conformal structure $X$ on the surface $S_g$, we consider Higgs bundles of the form $\mathcal{E}=\mathcal{K}^{1/2}\oplus \mathcal{O}\oplus \mathcal{K}^{-1/2}$ with Higgs field :
$$\Phi=\begin{pmatrix}
0 & t& 0 \\
0 & 0& t \\
\tau & 0 & 0 
\end{pmatrix}.$$

Here $\tau$ is the tautological section, and $t$ is any holomorphic section of $\mathcal{K}^{3/2}$. These Higgs bundles form a subset of the \emph{Slodowy slice} associated with the reducible representation $\iota_{\text{red}}:\mathfrak{sl}(2,\C)\to \mathfrak{sl}(3,\C)$, in the sense of \cite{magical}.

It is not known if these slices are disjoint if we vary the conformal structure $X$. If this was the case the union of the slices would form a subspace of dimension $10g-10$ of the character variety of expected dimension $16g-16$.

\medskip

Note that from the point of view of Higgs bundles, this slice of the character variety has a lot in common with the Hitchin component, which is the disjoint union for each Riemann surface structure $X$ of $S_g$ of a proper disk of dimension $10g-10$ in the character variety parametrized by a cubic differential. The fact that these slices are disjoint was conjectured and proven by Labourie \cite{Lab17}. 

Representations in the Hitchin component are the holonomies of convex projective structures, but also structures modeled on the space of full flags. These structures modeled on the space of full flags admit a natural fibration that can be defined using the associated affine sphere \cite{Loftin}.

\medskip

The non-Abelian Hodge correspondence allows to construct a representation from a Higgs bundle, but this construction involves solving a PDE, namely Hitchin's equation. It is difficult in general to read properties of the obtained representations on the Higgs bundles. For instance it is difficult to see which Higgs bundles correspond to discrete representations, or Anosov representations. However for this slice of Higgs bundles we obtain the following result:

\begin{thm}[{Theorem \ref{thm:mainText}}]
\label{thm:main}
Representations corresponding to the Higgs bundles $\mathcal{S}_X^\R$ are Anosov.
\end{thm}

Note that there is only one notion of Anosov representations in $\SL(3,\R)$, which are Borel-Anosov representations.

\medskip

A similar result was obtained in a work of Filip for compact slices of cyclic Higgs bundles, which are variations of Hodge structures, on a surface with or without boundary for the Lie group $\SO(2,3)$\cite{Fil21}. This result was extended to a larger non-compact slices containing these Hodge bundles by Zhang \cite{zhang}.

Our method involves understanding some geometric structure associated with the representation, which is also the case in \cite{Fil21}. Let $\mathcal{F}_{1,2}$ the space of full flags in $\R^3$. We prove in particular the following:

\begin{thm}[{Theorem \ref{thm:Fibration}}]
Representations corresponding to the Higgs bundles in $\mathcal{S}_X^\R$ are the holonomies of structures modeled on the space $\mathcal{F}_{1,2}$ of full flags in $\R^3$, on the unit tangent bundle of the surface, so that the fibers are conics in $\mathcal{F}_{1,2}$.
\end{thm}

A \emph{conic} in $\mathcal{F}_{1,2}$ is the set of points and tangent lines to a non-degenerate conic in $\mathbb{RP}^2$. Note that admitting such a structure is an open property so this does not characterize representations in $\mathcal{S}_X^\R$.

\medskip

We first show that the specific form of the Higgs bundle allows us to construct a geometric object, which is a parallel distribution of planes along the minimal immersion in the symmetric space $\mathbb{X}=\SL(3,\R)/\SO(3,\R)$. This is a bundle map $V:T\widetilde{S_g}\to T\mathbb{X}$ above a minimal immersion $\widetilde{S_g}\to \mathbb{X}$, naturally associated with the Fuchsian part $\Phi_0$ of the Higgs field, which is the part containing only the tautological section $\tau$ and not the differential $t$. For each unit vector $v\in T\widetilde{S_g}$ we define a open set $U_v\subset \mathcal{F}_{1,2}$, that we call a multicone, so that these multicones are nested along geodesics. This allows us to prove that the associated representation is Anosov thanks to a criterion we introduce in Theorem \ref{thm:CriterionAnosov}. The nestedness of these multicones gives indeed a control on the eigenvalue gap of elements in the representation.

\medskip

The multicones are constructed as follows: for every $x\in \widetilde{S_g}$ the plane $V(T_x\widetilde{S_g})$ is tangent to totally geodesic hyperbolic plane $\mathcal{H}_x$ in the symmetric space $\mathbb{X}$ that we call reducible, we define a map $p:\mathcal{F}_{1,2}\to \overline{\mathcal{H}_x}=\mathcal{H}_x\cup \partial \mathcal{H}_x$, which is a fibration over $\mathcal{H}$, that extends the nearest point fibration. The multicones we consider are pre-images by this map of half spaces in $\overline{\mathcal{H}_x}$ bounded by geodesics passing through $x$.

\medskip

The nestedness of the multicones, Theorem \ref{thm:Nestedness}, is the most technical part of the argument. The parallelism implies that this nestedness holds whenever $V(T_x\widetilde{S_g})$ is not too far from the tangent plane $\mathrm{d}h(T_x\widetilde{S_g})$, which amounts to controlling the norm of the holomorphic section $t$ for the harmonic metric. We prove such a bound on the norm of $t$ by a maximum principle, and check that this is exactly the control needed to obtain the nestedness property.

\medskip

The existence of this parallel distribution of planes describing the geometric structure is crucial in our argument. Parallel distributions of planes along a minimal surface in the symmetric space can be constructed similarly in rank $2$ for maximal representations and Hitchin representations but also for the non-maximal representations in $\SO(2,3)$ by Filip and Zhang \cite{Fil21,zhang}, as well as for some non-Hitchin $G_2'$-Higgs bundles studied by Collier and Toulisse \cite{collier2023holomorphiccurves6pseudospherecyclic}. 

\medskip

We begin this paper by introducing the slice in the character variety that we consider in Section \ref{sec:Slice}. We then introduce the Anosov property in Section \ref{sec:Anosov}, and prove a sufficient criterion to show this property using inclusions of multicones. We then show how to construct a parallel distribution of planes along the minimal surface in the symmetric space in Section \ref{sec:parallel}. In Section $5$ we study certain totally geodesic copies of $\mathbb{H}^2$ in the symmetric space, define a projection from the space of full flags onto these hyperbolic planes and define the multicones we want to consider. Finally in Section \ref{sec:Geom}, we show that the multicones associated to the parallel distribution of planes are nested, in order to prove the Anosov property and construct and a fibered geometric structure.

\subsection*{Acknowledgments}

We would like to thank Nicolas Tholozan for his helpful discussions on this subject at the beginning of the project.

S. Bronstein thanks the MPI-MIS for the quality of its working environment.

C. Davalo was funded by the RTG 2229 “Asymptotic
Invariants and Limits of Groups and Spaces” and the European Union via the ERC 101124349 "GENERATE". 
Views and opinions expressed are however those of the authors only and do not necessarily reflect those of the European Union or the European Research Council Executive Agency. Neither the European Union nor the granting authority can be held responsible for them.

%

\section{A slice of deformations of Barbot representations.}
\label{sec:Slice}

Before introducing the slice of the character variety that we consider in this paper, let us mention a few results on the non-Abelian Hodge correspondence.

\subsection{The non-Abelian Hodge correspondence for \texorpdfstring{$\SL(3,\R)$}{SL(3,R)}.}

Given a closed surface group $S_g$ of genus at least $2$, its \emph{character variety} $\chi(\pi_1(S_g), \SL(3,\R))$ is the geometric quotient of the space of representations:
 $$\text{Hom}(\pi_1(S_g),\SL(3,\R))// \SL(3,\R).$$
 In this expression, $\SL(3,\R)$ is acting by conjugation, and taking the geometric quotient means that we only consider semi-simple representations, or representations whose orbit by the action by conjugation is closed. In particular taking only this subset of the quotient makes it Haussdorf for the induced topology. 
 
 \medskip
 
The character variety can be parametrized by a moduli space of objects of analytic nature. This parametrization depends on the choice of a conformal structure $X$ on $S_g$, and is therefore not equivariant with respect to the mapping class group action. The objects in this moduli space are \emph{polystable $\SL(3,\R)$-Higgs bundles}, which are triples $(\mathcal{E},\Phi, B)$ where:

\begin{itemize}
\item $\mathcal{E}$ is a holomorphic rank $3$ vector bundle over $X$ with $\bigwedge^3 \mathcal{E}\simeq \mathcal{O}$,
\item $\Phi$ is a holomorphic section of $\mathcal{K}\otimes \text{End}(\mathcal{E})$ where $\mathcal{K}=\left(T_\C^{1,0}X\right)^*$ is the holomorphic cotangent line bundle of $X$,
\item $B$ is a holomorphic non-degenerate complex symmetric pairing of $\mathcal{E}$ so that $\Phi$ is $B$-symmetric.
\end{itemize}

Without the data of the bilinear pairing $B$, we have a $\SL(3,\C)$-Higgs bundle. 

\medskip

Such an object is called \emph{polystable} if $\mathcal{E}$ can be decomposed into a sum of holomorphic $\Phi$-invariant degree $0$ sub-bundles $(\mathcal{E}_i, \Phi_{|\mathcal{E}_i})$ that are \emph{stable}. 

A pair $(\mathcal{E}, \Phi)$ is \emph{stable} if every $\Phi$-invariant sub-bundle $\mathcal{F}$ of $\mathcal{E}$ has negative degree.

\medskip

Let $\mathcal{M}$ be the moduli space of polystable $\SL(3,\R)$-Higgs bundles up to  holomorphic gauge transformation.  

\begin{thm}[{\cite{Hit87},\cite{Cor88},\cite{Don03}}]
Given a polystable Higgs bundle $(\mathcal{E}, \Phi, B)$, there exists a metric $H$ on $\mathcal{E}$ whose Chern connection $\nabla$ satisfies Hitchin's equation:
 $$F_\nabla+[\Phi,\Phi^{*H}]=0.$$
 
The connection $D=\nabla+\Phi+\Phi^{*H}$ is flat and its holonomy defines an element of the character variety $\chi(\pi_1(S_g),\SL(3,\R))$. This defines a homeomorphism: $$\mathcal{M}\to \chi(\pi_1(S_g),\SL(3,\R)).$$ 
\end{thm}
\begin{rem}
	If the Higgs bundle $(\mathcal E, \Phi, B)$ is stable, then the metric is unique up to scale
	and its holonomy is irreducible in $\SL(3,\R)$. When $(\mathcal E, \Phi, B)$ is polystable,
	then the metric solution to Hitchin's equations decomposes along the decomposition
	into stable subbundles, and the metric is unique up to scale in each factor. In
	both cases, the holonomy of the metric does not depend on the choice of scaling involved.
\end{rem}

Therefore in order to construct representations up to conjugation, one can instead construct Higgs bundles.

\subsection{Deformations of Barbot representations.}

Let us fix a Riemann surface structure $X$ on the topological closed surface $S_g$ of genus $g$. Let us choose a square root $\mathcal{K}^{\frac{1}{2}}$ of the canonical bundle $\mathcal{K}=\left(T_\C^{(1,0)} X\right)^*$.

\medskip

These two choices determine a discrete and faithful representation $\rho:\pi_1(S_g)\to \SL(2,\R)$ up to conjugation, associated to the following $\SL(2,\R)$-Higgs bundle:

$$\mathcal{E}=\mathcal{K}^{\frac{1}{2}}\oplus\mathcal{K}^{-\frac{1}{2}}, \;\Phi=\begin{pmatrix}
0 & 0 \\
\tau & 0 
\end{pmatrix}.$$

Here $\tau\in  \Hom\left(\mathcal{K}^{\frac{1}{2}},\mathcal{K}^{-\frac{1}{2}}\right)\otimes \mathcal{K}\simeq \mathcal{O}$ is the tautological section, that is sometimes just denoted by $1$ in the literature.

The bilinear pairing on $\mathcal{E}$ that determines the real structure preserved by $\rho$ is :
$$B=\begin{pmatrix}
0 & 1 \\
1 & 0 
\end{pmatrix}.$$

\medskip

The corresponding $\SL(3,\C)$-Higgs bundle is the following 

$$\mathcal{E}=\mathcal{K}^{\frac{1}{2}}\oplus\mathcal{O} \oplus \mathcal{K}^{-\frac{1}{2}}, \;\Phi=\begin{pmatrix}
0 & 0& 0 \\
0 & 0& 0 \\
\tau & 0 & 0 
\end{pmatrix}.$$

This Higgs bundle is a fixed point of the $\mathbb{C}^*$-action. In particular it is a critical point of the energy functional. In order to deform this Higgs bundle so that its energy increases, one can consider deformations inside of the \emph{Slodowy slice} of the reducible $\mathfrak{sl}_2$ triple in $\mathfrak{sl}_3(\mathbb{R})$. Let $t_1,t_2$ a holomorphic sections of $\mathcal{K}^{\frac{3}{2}}$. We consider the following cyclic Higgs bundle:

$$\mathcal{E}=\mathcal{K}^{\frac{1}{2}}\oplus\mathcal{O} \oplus \mathcal{K}^{-\frac{1}{2}}, \;\Phi=\begin{pmatrix}
0 & t_1& 0 \\
0 & 0& t_2 \\
\tau & 0 & 0 
\end{pmatrix}.$$ 

We call $\mathcal{S}_X^\C$ the subset of the character variety associated to all of these Higgs bundles. 

\begin{prop}
Higgs bundles in $\mathcal{S}_X^\C$ are stable if $t_1,t_2\neq 0$ and polystable if $t_1,t_2=0$.
\end{prop}

\begin{proof}
	Assume $t_1,t_2\neq 0$
	Let $L$ be a $\Phi$-invariant line subbundle of $\calE$ and
	consider a meromorphic section $(a,b,c)$ of $\calE$.
	Then there must be a holomorphic section $\lambda$ of $\calK$,
	such that:
	\begin{equation}
		\lambda(a,b,c)=(t_1b,t_2c,a)\,.
	\end{equation}
	But then, either $b$ identically vanishes, or $\lambda^3=s^2$.

	If $b=0$, then $s c=0$, so $c=0$ and $a=0$ too, which is absurd.
	In the other case, generically there is no solution to $\lambda^3=t_1t_2$, hence there is no invariant subbundle.
	If there is a $\lambda$ solution, then an appropriate section of
	the corresponding line bundle would be $(\lambda^2/t_1,1,\lambda/t_2)$.
	Since $\lambda/t_1$ is a meromorphic section of $\calK^{-\frac{1}{2}}$,
	we see that our meromorphic section has no zero but at least
	$g-1$ poles, hence the invariant line bundle is of negative degree,
	as desired.

	When $t_1,t_2=0$, the bundle decomposes as the sum $\calK^{\frac{1}{2}}\oplus\calK^{-\frac{1}{2}}$, and $\calO$. Both are stable degree $0$ line bundles,
	so it is polystable, as desired.
\end{proof}

Note that for every $\lambda\in \C$ there is a Gauge transformation of this Higgs bundle that maps the Higgs field $\Phi$ to a similar Higgs field with $(t_1,t_2)$ replaced by $(\lambda t_1,\lambda^{-1}t_2)$. Conversely two such Higgs bundles are gauge equivalent if and only if they are of this form, since $\tr(\Phi^3)=t_1t_2$ is a Gauge invariant.

\begin{que}
Are the slices $\mathcal{S}_{X_1}^\C$ and $\mathcal{S}_{X_2}^\C$ disjoint for all distinct pairs of Riemann surface structures $X_1,X_2$?
\end{que}

When $t_1=t_2=t$, these Higgs bundles preserve the bilinear pairing :
$$B=\begin{pmatrix}
0 & 0& 1 \\
0 & 1& 0 \\
1 & 0& 0 
\end{pmatrix}.$$

This implies that they define representations $\rho:\pi_1(S_g)\to \SL(3,\R)$. We denote by $\mathcal{S}^\R_X$ the corresponding subset of the character variety of $\SL(3,\R)$.

\begin{prop}
The slice $\mathcal{S}_X^\R$ is topologically a properly embedded disk of real dimension $4g-4$
\end{prop}

\begin{proof}
Since $\mathcal{K}\otimes \mathcal{K}^{-\frac{3}{2}}$ has negative degree, the complex dimension
of $H^0(\mathcal{K}^{3/2})$ is equal to, by the Riemann--Roch formula:
$$\deg(\mathcal{K}^{3/2})+1-g=2g-2.$$

The application that associates to a section of $\mathcal{K}^{\frac{3}{2}}$ up to sign associates the Higgs bundle $(\mathcal{E},\Phi)$ in the moduli space is injective and proper, since the invariant polynomial $\text{tr}(\Phi^2)$ is equal to $t^2$. Indeed, for $t'\neq t,-t$, the Higgs bundles cannot be gauge equivalent, and if $t$ goes to $+\infty$ they diverge in the character variety.

It remains to check that the $\SL(3,\R)$-Higgs bundles are not Gauge--conjugate for $t$ and $-t$.
	This is a direct computation that we do in Proposition~\ref{prop:Gauge}
\end{proof}
Let us now compare this slice with the deformations considered by Barbot. Let us consider the following general $\SL(3,\R)$-Higgs bundle $(\mathcal{E},\Phi_{t,\delta,q})$:

$$\mathcal{E}=\mathcal{K}^{\frac{1}{2}}\oplus\mathcal{O} \oplus \mathcal{K}^{-\frac{1}{2}}, \;\Phi=\begin{pmatrix}
\delta & t& q \\
0 & -2\delta& t \\
\tau & 0 & \delta 
\end{pmatrix}.$$

Here $\tau$ is the tautological section, $\delta$, $t$ and $q$ are respectively holomorphic sections of $\mathcal{K}$, $\mathcal{K}^{\frac{1}{2}}$ and $\mathcal{K}^2$. These Higgs bundles form the intersection of the \emph{Slodowy slice} associated with the reducible $\mathfrak{sl}_2$ triple in $\mathfrak{sl}_3(\C)$, following the terminology of \cite{magical}, with the space of $\SL(3,\R)$-Higgs bundles. 

\medskip

Fuchsian representations, for a given choice of a square root $\mathcal{K}^{\frac{1}{2}}$, correspond to such Higgs bundles with $t=\delta=0$. The deformations considered by Barbot \cite{Bar10} correspond exactly to the case $t=0$. Note that the choice of $\delta$ corresponds to choosing an element in $H^0(\mathcal{K})\simeq H^1(\mathcal{O})\simeq H^1(\Gamma_g, \R)$. 

\medskip

The following Proposition ensures that our slices are embedded in the character variety, and that this slice intersects the slice considered by Barbot only at the Fuchsian locus.

\begin{prop}
\label{prop:Gauge}
If a $\SO(3,\C)$ holomorphic gauge transformation conjugates the Higgs bundle $(\mathcal{E},\Phi_{t,\delta,q})$ into $(\mathcal{E},\Phi_{t',\delta',q'})$, then it is diagonal and $t'=t$, $\delta'=\delta$, $q'=q$.
\end{prop}
\begin{proof}
Let $g$ be such a holomorphic gauge transformation of $\mathcal{E}$ conjugating $\Phi_{t,\delta,q}$ into $\Phi_{t',\delta',q'}$.  It has to be upper diagonal because its lower diagonal coefficients define holomorphic sections of $\mathcal{K}^{-\frac{1}{2}}$ and $\mathcal{K}^{-1}$, which are line bundles of negative degree. Let us write the gauge transformation for some sections $a,b,c,d,e,f$ as :

$$g=\begin{pmatrix}
a & e& d \\
0 & b& f \\
0 & 0 & c
\end{pmatrix}.$$

If this gauge transformation conjugates the two Higgs fields, it must satisfy the following equation:
$$g\Phi_{t,\delta,q}=\begin{pmatrix}
* & *& * \\
\tau f & *& * \\
\tau c & 0 & *
\end{pmatrix} = \begin{pmatrix}
* & *& * \\
0 & *& * \\
\tau c & \tau e & *
\end{pmatrix}=\Phi_{t',\delta',q'}g.$$


In particular $e=f=0$. Finally the fact that $g$ lies in $\SO(3,\C)$, \ie preserve the pairing $B$, implies that it is diagonal of the form $\text{Diag}(\lambda, 1,\lambda^{-1})$. Necessarily $\lambda=1$, hence $t'=t$, $\delta'=\delta$, $q'=q$. 
\end{proof}

\begin{rem}

Note here that the gauge transformation of $\mathcal{E}$ given by $\text{Diag}(\omega,\omega^{-2},\omega)$ which is in $\SL(3,\C) \setminus \SO(3,\C)$ where $\omega^6=1$ is a primitive root of unit that conjugates $\Phi_{t,\delta,q}$ into $\Phi_{-t,\delta, q}$. Therefore the representations associated with $(\mathcal{E},\Phi_{t,\delta,q})$ and $(\mathcal{E},\Phi_{-t,\delta,q})$ are conjugated in $\SL(3,\C)$. However these representations are not conjugated in $\SL(3,\R)$ because of Proposition \ref{prop:Gauge}. 

\end{rem}

\section{The Anosov property.}
\label{sec:Anosov}

Anosov representations, which were introduced by Labourie \cite{Labourie} and defined in more generality by  Guichard-Wienhard \cite{GW12}, are representations of hyperbolic groups  that are discrete and faithful in a strong sense. For representations into $\SL(3,\R)$ there is only one notion of Anosov representations, which are Borel Anosov representations.

\medskip

We will use the following definition, due to Bochi-Potrie-Sambarino \cite{BPS} and Kapovich-Leeb-Porti \cite{KLP17}. Let $\Gamma$ be a Gromov hyperbolic group. Let us fix a scalar product on $\R^3$, and let $\sigma_i(g)$ be the $i$-th largest modulus of the sigular values of $g\in \SL(3,\R)$. Let us fix a word metric $|\cdot|$ on $\Gamma$. This definition, is independent of these choices.

\begin{defn}
A representation $\rho:\Gamma\to SL(3,\R)$ is Anosov if for some $A,B>0$, for all $\gamma\in  \Gamma$:
$$\sigma_1(\rho(\gamma))-\sigma_2(\rho(\gamma))\geq A |\gamma|-B.$$
\end{defn}

Here $|\gamma|$ is the word lenth of $\gamma$ for any choice of generators of $\Gamma$.
\begin{rem}
Note that $\sigma_1(g)-\sigma_2(g)=\sigma_2(g^{-1})-\sigma_3(g^{-1})$, so Anosov representations also satisfy $\sigma_2(\rho(\gamma))-\sigma_3(\rho(\gamma))\geq A |\gamma|-B$.
\end{rem}

Anosov representations $\rho:\Gamma\to\SL(3,\R)$ come with a unique continuous and $\rho$ equivariant map \emph{limit map} $\xi_\rho=(\xi^1_\rho,\xi^2_\rho):\partial \Gamma\to \mathcal{F}_\Delta$ into the space of full flags in $\R^3$ where $\partial \Gamma$ is the Gromov bondary of $\Gamma$, and the image of the attracting fixed point $\gamma^+\in \partial \Gamma$ of any hyperbolic element $\gamma\in \Gamma$ is mapped by $\xi_\rho$ to the attracting fixed flag of $\rho(\gamma)$.

\medskip

In the case of a surface group, an equivalent characterization is the following. Let $\lambda_i(g)$ be the logarithm of the $i$-th largest characteristic value of $g$ (i.e eigenvalue if $g$ is diagonalizable). 

\begin{thm}[{\cite{Kostas}}]
\label{thm:CharacAnosovEigen}

A representation $\rho:\Gamma_g\to SL(3,\R)$ is Anosov if for some $A,B>0$, for all conjugacy classes $\overline{\gamma}$ in $\Gamma$:
$$\lambda_1(\rho(\gamma))-\lambda_2(\rho(\gamma))\geq A \ell(\overline{\gamma})-B.$$
\end{thm}

Here $\ell(\overline{\gamma})$ is the length of the smallest element in the conjugacy class for the chosen word metric on $\Gamma$ . If $\Gamma=\pi_1(S_g)$ one can replace it with the length of the shortest geodesic representative of the associated simple closed curve, for a fixed Riemannian metric on $S_g$.

\subsection{Multicone criterion}
\label{subsec:Multicone}
We will use a criterion to prove the Anosov property for representations of surface groups that uses multicones, similar to a criterion introduced Bochi-Potrie-Sambarino \cite{BPS}. This criterion can be seen as a continuous, but compact, analog of a ping-pong lemma.

\medskip

Let $\mathcal{F}_{1,2}=\mathcal{F}_\Delta$ be the set of full flags in $\R^3$, \ie pairs $f=(f^1,f^2)$ where $f^1\subsetneq f^2 \subsetneq R^3$. We write $\mathcal{F}_{1,2}$ for the space in which the multicones will live and $\mathcal{F}_\Delta$ for the space in which the limit map lives: this distinction is useful for clarity, as well as for considering generalizations of the criterion. 

To $f=(f^1,f^2)\in \mathcal{F}_\Delta$ we associate its \emph{thickening} $K_f\subset \mathcal{F}_{1,2}$ as the set of flags $g=(g^1,g^2)\in \mathcal{F}_{1,2}$ such that $f^1=g^1$ or $g^2=f^2$. These thickenings have the following basic property:

\begin{prop}
Two flags $f_1,f_2\in \mathcal{F}_\Delta$ are transverse if and only if $K_{f_1}\cap K_{f_2}=\emptyset$ and are equal if and only if $K_{f_1}=K_{f_2}$.
\end{prop}

We say that an open subset $U\subset \mathcal{F}_{1,2}$ is a \emph{multicone} if there exist $f_+, f_-\in \mathcal{F}_\Delta$ such that $K_{f_+}\subset U$ and $K_{f_-}\cap \overline{U}=\emptyset$. 

\begin{rem}
For a more general semi-simple Lie group $G$, given two flag manifolds $\mathcal{F}_\Theta$, $\mathcal{F}$ playing the roles of $\mathcal{F}_\Delta$ and $\mathcal{F}_{1,2}$, one can define a notion of multicones in $\mathcal{F}$ by defining a notion of thickening $K_f\subset \mathcal{F}$ for $f\in \mathcal{F}_\Theta$. For that one just need to have a balanced Tits-Bruhat ideal, see \cite{KLP}. Our criterion naturally extends to a criterion that certifies the $\Theta$-Anosov property.
\end{rem}

In order to topologize the space of multicones, we define $d(U_1,U_2)=d_H(\overline{U_1},\overline{U_2})+d_H(\mathcal{F}_{1,2}\setminus U_1,\mathcal{F}_{1,2}\setminus U_2)$ where $d_H$ is the Haussdorf distance for some Riemannian metric on $\mathcal{F}_{1,2}$.
In particular for a fixed closed subset $K\subset \mathcal{F}_{1,2}$ , the sets $\lbrace U \text{ multicone}\mid K\subset U\rbrace$ and $\lbrace U \text{ multicone}\mid K\cap \overline{ U}=\emptyset \rbrace$ are open.

\medskip

We can show that a representation is Anosov by presenting a compact collection of multicones with certain nestedness properties. 

\begin{thm}
\label{thm:CriterionAnosov}
Let $\rho:\Gamma_g\to \SL(3,\R)$ be a representation. Fix a Riemannian metric on $S_g$. Suppose that one can associate continuously and $\rho$-equivariantly to every unit tangent vector $v\in T^1\widetilde{S_g}$ in the universal cover a multicone $U_v\subset \mathcal{F}_{1,2}$ so that for all $t>0$, $\overline{U_{\Phi_t(v)}}\subset U_v$ where $\Phi$ is the geodesic flow on $T^1\widetilde{S_g}$. Then the representation $\rho$ is Borel-Anosov. 

Moreover the limit map $\xi_\rho:\partial \Gamma\to \mathcal{F}_\Delta$ is characterized by the fact that for all $v\in T^1\widetilde{S_g}$ tangent to a closed curve, if $x=\lim_{t\to \infty} \Phi_t(v)$:
 $$K_{\xi_\rho(x)}=\bigcap_{t>0} U_{\Phi_t(v)}.$$
 
 If the metric on $S_g$ is negatively curved, this holds for all $v\in T^1\widetilde{S_g}$.
\end{thm}

\begin{rem}
Given a multicone $U_v\subset \mathcal{F}_{1,2}$ one can define a non-empty open set $V_v=\lbrace f\mid K_f\subset U_v\rbrace \subset \mathcal{F}_\Delta$. These sets will have the same nestedness properties, and Theorem \ref{thm:CriterionAnosov} is stating that the intersection of the $V_{\Phi_t(v)}$ is equal to $\lbrace \xi_\rho(x)\rbrace$.
\end{rem}

In order to prove that the representation is Anosov we construct its boundary map. For that we need to measure how much two multicones are nested into one another.

\begin{defn}
Let $U_1, U_2 \subset \mathcal{F}_{1,2}$ be two multicones with $\overline{U_2}\subset U_1$. The \emph{nestedness} $\text{Nest}(U_1,U_2)$ is the minimum for all flags $f_+,f_-\in \mathcal{F}_\Delta$ such that $K_{f_+}\subset \overline{U_2}$ and $K_{f_-}\cap U_1=\emptyset$ of the largest $\lambda\in \R^{>0}$ such that $U_2\subset g_\lambda\cdot U_1$ where $g_\lambda\in \SL(3,\R)$ can be written as follows in some basis associated to $f_+,f_-$:

\begin{equation}
\label{eq:g lambda}
g_\lambda=\begin{pmatrix}
e^\lambda & 0& 0 \\
0 & 1& 0 \\
0 & 0& e^{-\lambda} 
\end{pmatrix}.
\end{equation}
\end{defn}

Note that for a fixed $f_+,f_-$ and $\lambda$, there is a single such $g_\lambda$. Moreover the set of $f^+$ such that $K_{f^+}\subset \overline{U_2}$ is compact. In particular if $\overline{U_2}\subset U_1$ for multicones $U_1,U_2$, then $\text{Nest}(U_1,U_2)>0$.

This quantity satisfies the following inequality.

\begin{prop}
For every multicones $U_1,U_2,U_3\subset \mathcal{F}_{1,2}$ such that $\overline{U_3}\subset U_2$ and $\overline{U_2}\subset U_1$, $\text{Nest}(U_1,U_3)\geq \text{Nest}(U_1,U_2)+ \text{Nest}(U_2,U_3)$.
\end{prop}
\begin{proof}
	Let $(f_+,f_-)$ be two flags such that $K_{f_+}\subset\overline{U_3}$
	and $K_{f_-}\cap U_1=\emptyset$.
	Let then $\lambda_{ij}$ be the maximal value such
	that $U_j\subset g_{\lambda_{ij}}U_i$ for all $g_{\lambda_{ij}}$ which can be written as previously in some basis adapted to $f_+$ and $f_-$. Because $g_{\lambda_{13}}$ can be taken to be equal to $g_{\lambda_{23}}g_{\lambda_{12}}$, and since $U_3\subset g_{\lambda_{23}}g_{\lambda_{12}}U_1$:
	\begin{equation*}
		\lambda_{13}\geq\lambda_{12}+\lambda_{23}\,.
	\end{equation*}
	By superadditivity of the minimum, we get
	\begin{equation*}
		\text{Nest}(U_1,U_3)\geq\min\lambda_{12}+\min \lambda_{23}\,.
	\end{equation*}
	
	In this expression the minima are taken over $\lbrace f_+,f_-\mid K_{f_+}\subset \overline{U_3},K_{f_-}\cap U_1=\emptyset\rbrace$.
	But this set over which this minimum is taken is strictly smaller than
	the ones in the definition of the nestedness, hence the last inequality:
	\begin{equation*}
		\text{Nest}(U_1,U_3)\geq \text{Nest}(U_1,U_2)+\text{Nest}(U_2,U_3)\,.
	\end{equation*}
\end{proof}

\begin{prop}
\label{IntersectionMulticones}
Let $(U_n)$ be a sequence of multicones in $\mathcal{F}_{1,2}$ such that for all $n\in \mathbb{N}$, $\overline{U_{n+1}}\subset U_n$. If $\text{Nest}(U_1,U_n)$ tends to $+\infty$ when $n$ goes to $+\infty$, there exist a single flag $f\in \mathcal{F}_\Delta$ such that:
$$K_f=\bigcap_{n\in \mathbb{N}}U_n .$$
\end{prop}

\begin{proof}
For each $n\in \mathbb{N}$ let $\lambda_n=\text{Nest}(U_1,U_n)$ and $g_n\in \SL(3,\R)$ be the element $g_{\lambda_n}$ from the definition of the nestedness distance. The attracting and repulsive fixed flags of $g_n$ converge up to subsequence to some flag $f_+,f_-$ such that $K_{f_-}\cap \overline{U_1}=\emptyset$ and :
$$K_{f_+}\subset\bigcap_{n\in \mathbb{N}}U_n.$$

In particular the flags $f_+$ and $f_-$ are transverse. Since $\lambda_n$ converges to $+\infty$, for every element in $p\in\mathcal{F}_{1,2}\setminus K_{f_-}$ the sequence $(g_n\cdot p)$ converges up to subsequence to a point in $K_{f_+}$. Hence $K_{f_+}\subset\bigcap_{n\in \mathbb{N}}U_n$. This element is unique since a flag is characterized by its thickening.
\end{proof}

An application of this is the following lemma :

\begin{lem}
\label{lem:Eigenvalue}
Let $g\in G$ be an element with a fixed flag $f\in \mathcal{F}_\Delta$. Suppose that for some multicone $U\subset \mathcal{F}_{1,2}$, $K_f\subset U$ and $g\cdot \overline{U}\subset U$. In a basis such that matrices fixing $f$ are upper triangular, if the diagonal coefficients of $g$ are $(\pm e^{\lambda_1},\pm e^{\lambda_2}, \pm e^{\lambda_3})$, then $\lambda_1-\lambda_2\geq 0$ and $\lambda_2-\lambda_3\geq 0$.
\end{lem}

\begin{proof}
If $\lambda_1-\lambda_2< 0$ or $\lambda_2-\lambda_3< 0$ then $f$ is not an attracting fixed point of the action of $g$ on $\mathcal{F}_\Delta$. We therefore show that $g$ has $f$ as an attracting fixed point. The nestedness $\text{Nest}(U, g^n U)$ is at least equal to $n\times \text{Nest}(U, g U)>0$ hence it goes to $+\infty$. in particular the intersection of $g^n U$ for all $n$ is a singleton, necessarily equal to $f$. Hence we have proven that $f$ is an attracting fixed point of $g$.
\end{proof}

Now we can prove our main characterization, using Theorem \ref{thm:CharacAnosovEigen}.
\begin{proof}[{Proof of Theorem \ref{thm:CriterionAnosov}}]
Let $\delta>0$ be the infimum of $\text{Nest}(U_v, U_{\Phi_1(v)})$ for $v\in T^1\widetilde{S_g}$. We know that $\text{Nest}(U_v, U_{\Phi_t(v)})> \delta (t-1)$ for all $t>0$. In particular we can apply Proposition \ref{IntersectionMulticones} and define for all $v$ a pair of transverse flags $f^+(v),f^-(v)$ such that 
$$K_{f^+(v)}=\bigcap_{t>0}U_{\Phi_t(v)}, $$
$$K_{f^-(v)}=\bigcap_{t<0}\mathcal{F}_{1,2}\setminus \overline{U_{\Phi_{t(v)}}}.$$

\medskip

Let $\overline{\gamma}$ be a simple closed curve, i.e. a conjugacy class in $\pi_1(S_g)$, with shortest geodesic representative $\eta$ of length $\ell$. Let $v\in T^1\widetilde{S_g}$. Let $\gamma$ be the representative of $\overline{\gamma}$ sending $v$ to $\Phi_{\ell}(v)$. Since $\rho(\gamma^n)\overline{U_v}=U_{\Phi_{n\ell}(v)}\subset U_v$, we know that $\rho(\gamma)$ fixes $f_+(\gamma)$.

Let $g=g_{\lambda}$ be as in \eqref{eq:g lambda} for the flags $f_+(v)$ and $f_-(v)$ and for $\lambda=\delta(\ell-1)$. Since $\text{Nest}(U_v, U_{\Phi_1(v)})>\lambda$, then one an take $g$ so that $U_{\Phi_1(v)}\subset gU_v$. The product $g_t^{-1}\rho(\gamma)$ fixes $f^+(v)$ and satisfies $g_\lambda^{-1}\rho(\gamma)U_v\subset U_v$. Therefore by Lemma \ref{lem:Eigenvalue} the diagonal coefficients $(\pm e^{\lambda_1},\pm e^{\lambda_2}, \pm e^{\lambda_3})$ of the upper-diagonal matrix in the basis adapted to $f^+(v)$ corresponding to $\rho(\gamma)$ satisfy:
$\lambda_1-\lambda_2-\lambda\geq 0$, $\lambda_2-\lambda_3-\lambda\geq 0$. In particular $\lambda_1-\lambda_2\geq \delta\ell-\delta$ so $\rho$ is Anosov thanks to Theorem \ref{thm:CharacAnosovEigen}.

Moreover $f^+(\gamma)$ is the attracting fixed point of $\gamma$, so it is equal to $\xi_\rho(x)$ for $x$ the limit of $(\Phi_t(v))$.

\medskip

Finally for the last point, note that if the metric is negatively curved then the set of vectors tangent to closed curves is dense. Hence by the continuity of the limit map $f^+(v)=\xi_\rho(x)$.
\end{proof}

\section{A parallel distribution of planes.}
\label{sec:parallel}
In this section we use the information we have on our special Higgs bundles to obtain a geometric object in the symmetric space: a parallel distribution of planes along a minimal surface in the Riemannian symmetric space $\mathbb{X}=\SL(3,\R)/\SO(3,\R)$.
\subsection{The Higgs bundle, and the parallel distribution of planes.}
We consider the Higgs bundle $\mathcal{E}=\mathcal{K}^{1/2}\oplus \mathcal{O}\oplus\mathcal{K}^{-1/2}$ on the surface $S_g$ with Riemann surface structure $X$ with Higgs field:
$$\Phi= \begin{pmatrix}
0 & t & 0\\
0 & 0 & t\\
\tau & 0 & 0
\end{pmatrix}.$$ 

The real form is determined by the following complex bilinear pairing:

$$B=\begin{pmatrix}
0 & 0 & 1\\
0 & 1 & 0\\
1 & 0 & 0
\end{pmatrix}.$$

The harmonic metric, characterized by Hitchin's equation, is diagonal since the Higgs field is cyclic \cite{sim06}:
$$H=\begin{pmatrix}
h_1 & 0 & 0\\
0 & 1 & 0\\
0 & 0 & h_1^{-1}
\end{pmatrix}.$$ 

The Riemannian symmetric space $\mathbb{X}$ associated with $\SL(3,\R)$ can be seen as the space $\SL(3,\R)/\SO(3,\R)$ of scalar products of volume one on $\R^3$. The harmonic metric, together with the associated flat connection define a map $h:\widetilde{S_g}\to \mathbb{X}$ that is equivariant with respect to the holonomy representation of the flat connection. The differential of this map in a direction $v\in T\widetilde{S_g}$ corresponds via the Maurer-Cartan form to $\Psi(v)=\Phi\left(v^{(1,0)}\right)+\Phi^*\left(v^{(1,0)}\right)$. Via this identification, if $D$ is the trivial flat connection on $\mathfrak{sl}(3,\R)\times\widetilde{S}$ and $\nabla^\mathbb{X}$ the Levi-Civita connection on $T\mathbb{X}$ which is identified as a subbundle of $\mathfrak{sl}(3,\R)\times \widetilde{S_g}$ by the Mauer-Cartan form:
$$D=\nabla^\mathbb{X} +\text{ad}_{\Psi(v)}.$$
In particular $\nabla^\mathbb{X}$ coincides with the Chern connection $\nabla^H$ associated with the harmonic map $h$ on the flat bundle $\mathfrak{sl}(3,\R)\times\widetilde{S_g}$ identified as a subbundle of $\text{End}(\mathcal{E})$.

 \medskip

We consider the following other section of $\text{End}(\mathcal{E})\otimes \mathcal{K}$:

$$\Phi_0= \begin{pmatrix}
0 & 0 & 0\\
0 & 0 & 0\\
\tau & 0 & 0
\end{pmatrix}.$$ 

The data of $\Phi_0+\Phi_0^*$ corresponds via the Maurer-Cartan map to a vector bundle homomorphism $V:T\widetilde{S_g}\to T\mathbb{X}$ over $h:\widetilde{S_g}\to \mathbb{X}$, but that is not the tangent map.
The image of $V$ defines a distribution of tangent planes along the harmonic map $h:\widetilde{S_g}\to X$. 

\medskip

The harmonic map $h$ defines an induced metric on $X=S_g$, but the metric that will play the main role for us is an other metric on $X$ defined in a natural way using $H$. The metric $H$ defines a metric $h_1^{-2}$ on the line bundle $\left(\mathcal{K}^{1/2}\right)^*\otimes \mathcal{K}^{-1/2}\simeq \mathcal{K}^{-1}$. this defines a metric on $X$, that is the one we consider from now on.

\medskip

The following proposition describes the crucial property of this distribution of planes.

\begin{prop}
\label{prop:ParallelPlanes}
The distribution of planes $V(T\widetilde{S_g})$ along $h:\widetilde{S_g}\to \mathbb{X}$ is parallel for the Levi-Civita conection on $\mathbb{X}$. More precisely one has for every section $v$ of $T\widetilde{S_g}$:
$$\nabla^\mathbb{X} V(v) =V(\nabla^{X} v).$$

In this expression $\nabla^X$ is the Levi-Civita connection for the metric $h_1^{-2}$ on $X=S_g$, and $\nabla^{\mathbb{X}}$ is the Levi-Civita connection for the Riemannian metric on $\mathbb{X}$.
\end{prop}

As a consequence if $v$ is a parallel section of $T\widetilde{S_g}$ over a curve $\gamma$ for the connection $\nabla^X$, then $V(v)$ is a parallel section of $T\mathbb{X}$ over the curve $h\circ\gamma$.

\begin{rem}
The Higgs field satisfies an equation that is often written as $\nabla^H \Phi=0$, which should be understood as a (skew-symmetric) $2$-form. This notation can be slightly misleading here as this does not imply for all sections $v$ of $T\widetilde{S_g}$ that $\nabla^H (\Phi\left(v^{(1,0)}\right))=\Phi\left(\nabla v^{(1,0)}\right)$ viewed as sections of $\text{End}(\mathcal{E})$.
\end{rem}

\begin{proof}
By definition $V(v)$ can be identified with $\Phi_0\left(v^{(1,0)}\right)+\Phi_0^{*}\left(v^{(1,0)}\right)$ where $v^{(1,0)}$ is the $(1,0)$-part of $v$. The element $\nabla^\mathbb{X} V(v)\in T\mathbb{X}$ can be identified with the following $H$-hermitian endomorphism of $\mathcal{E}$:
$$ \nabla^H \left(\Phi_0\left(v^{(1,0)}\right)+\Phi_0^{*}\left(v^{(1,0)}\right)\right)=\nabla^H\begin{pmatrix}
0 & 0 & \tau^*\left(v^{(1,0)}\right)\\
0 & 0 & 0\\
\tau\left(v^{(1,0)}\right) & 0 & 0
\end{pmatrix}.$$

The line decomposition of $\mathcal{E}$ is preserved by $\nabla^H$ since the harmonic metric $H$ is diagonal. Note also that the metric induced by $H$ on $\left(\mathcal{K}^{1/2}\right)^*\otimes \mathcal{K}^{-1/2}\subset \End(\mathcal{E})$ is $h_1^{-2}$ and that $\tau$ is the tautological section. Therefore $\nabla^\mathbb{X} V(v)$ corresponds to:

$$\begin{pmatrix}
0 & 0 & \nabla^H \tau^*\left(v^{(1,0)}\right)\\
0 & 0 & 0\\
\nabla^H \tau\left(v^{(1,0)}\right) & 0 & 0
\end{pmatrix}=\begin{pmatrix}
0 & 0 & \tau \left(\left(\nabla^Xv\right)^{(1,0)}\right)^*\\
0 & 0 & 0\\
\tau \left(\left(\nabla^Xv\right)^{(1,0)}\right) & 0 & 0
\end{pmatrix}$$

This last term is equal to $\Phi_0\left(\left(\nabla^Xv\right)^{(1,0)}\right)+\Phi_0^{*}\left(\left(\nabla^Xv\right)^{(1,0)}\right)$, which is in turn identified with $V(\nabla^{X} v)$
\end{proof}

Note that this distribution of planes coincides with the tangent distribution of planes to the minimal surface exactly at the zeroes of the section $t\in H^0(\mathcal{K}^{\frac{3}{2}})$.

%
%
%

\subsection{The maximum principle}

We apply a maximum principle to show that the distribution of planes $V$ is not too far from the tangent distribution of planes of the minimal surface. Concretely we control the norm of the section $t$ by of the norm of the tautological section $\tau$.

\begin{prop}
\label{prop:Maximal Principle}
Let $\lVert \tau \rVert$ and $\lVert t \rVert$ be the norms of the sections $\tau$ and $t$ with respect to the harmonic metric $H$ and any metric on $X$. One has $\lVert t \rVert<\lVert \tau \rVert.$
\end{prop}

\begin{proof}
	Denote by $\omega$ the Volume form of the chosen metric on $X$.
	Since $t$ is a $(1,0)$-form, one has the identification $t\wedge t^\ast=-2i\|t\|^2\omega$.

Note that the metric on $X$ does not influence the ratio $\lVert t \rVert/\lVert \tau \rVert$. 
Consider $(\calE,\dbar,\Phi)$ the previously considered Higgs bundle. Hitchin's equation write
as:
\begin{equation}
\label{eq:Hitchin}
	\left\{\begin{array}{ll}
	F_{\calK^\frac{1}{2}}+t\wedge t^\ast-\tau\wedge\tau^\ast&=0\\
	F_\calO&=0\,.
	\end{array}\right.
\end{equation}
Since $t$ and $\tau$ are holomorphic sections of respectively $\left(\mathcal{K}^{-\frac{1}{2}}\right)^*\otimes \mathcal{O}\otimes \mathcal{K}$ and $\mathcal{K}^2= \left(\mathcal{K}^{\frac{1}{2}}\right)^*\otimes\mathcal{K}^{-\frac{1}{2}}\otimes \mathcal{K}$, if we compute their norms for the harmonic metric $H$ and any metric on $\mathcal{K}$:
\begin{equation*}
	\left\{\begin{array}{rl}
	\partial\dbar \log(\lVert t \rVert)&=
	F_{\mathcal{K}^{-\frac{1}{2}}}-F_{\mathcal{K}}=t\wedge t^\ast-\tau\wedge\tau^\ast-F_{\mathcal{K}}\\
	\partial\dbar \log(\lVert \tau \rVert)&=
	F_{\mathcal{K}^{-\frac{1}{2}}}-F_{\mathcal{K}^{\frac{1}{2}}}-F_{\mathcal{K}}=-2t\wedge t^\ast+2\tau\wedge\tau^\ast-F_{\mathcal{K}}\,.
	\end{array}\right.
\end{equation*}

Now we consider the difference of these, that does not depend on the metric on $\mathcal{K}$:
$$\partial\dbar\log(\lVert t \rVert/\lVert \tau \rVert)=3t\wedge t^\ast-3\tau\wedge\tau^\ast.$$

Let $\widetilde\Delta=\frac{\partial\wedge\dbar}{\tau\wedge\tau^\ast}$ be the Laplacian
with regard to the metric $h_1^{-2}$ on $X$. For that metric, we have the following elliptic PDE:
\begin{equation*}
	\widetilde\Delta \log\frac{\|t\|}{\|\tau\|}=3(\frac{\|t\|^2}{\|\tau\|^2}-1)\,.
\end{equation*}
Applying the strong maximum principle to that PDE, we obtain that at a maximum of $\|t\|/\|\tau\|$,
$\|t\|\leq\|\tau\|$ with equality if and only if $\|t\|=\|\tau\|$ everywhere, which never happens
since $t$ must vanish at some point.
We deduce that $\|t\|<\|\tau\|$, as desired.

\end{proof}

\begin{rem}
\label{rem:NegativeCurved}

Note that a consequence of the maximum principle and Hitchin's equation \eqref{eq:Hitchin} is that the metric $h_1^{-2}$ is negatively curved on $X$
\end{rem}

\section{Construction of the multicones}
\label{sec:MulticonesConstruction}

In this section we construct a multicone associated to a reducible hyperbolic plane $\mathcal{H}$ in the symmetric space and a tangent vector. For that we define a projection of the flag manifold $\mathcal{F}_{1,2}$ onto the compactification of $\mathcal{H}$ in $\mathbb{X}$. We then describe in more detail this fibration, and define the multicone as the pre-image of a half plane of $\mathcal{H}$. 


\subsection{Fibration over a reducible hyperbolic plane}
We say that a totally geodesic complete surface $\mathcal{H}$ in the symmetric space $\mathbb{X}$ associated to $\SL(3,\R)$ is a \emph{reducible $\mathbb{H}^2$} if it is a copy of the hyperbolic plane equivariant with respect to a reducible and non-trivial representation of $\mathfrak{sl}(2,\R)\to \mathfrak{sl}(3,\R)$. In other words it is equal up to the action of an element of $\SL(3,\R)$ to:
$$\mathcal{H}_0=\left\lbrace \begin{pmatrix}
a & 0 & c\\
0 & 1 & 0\\
c & 0 & b
\end{pmatrix}\mid ab-c^2=1\right\rbrace.$$

Note that here we identify $\mathbb{X}=\SL(3,\R)/\SO(3,\R)$ with the space of scalar products of volume $1$, or symmetric matrices of determinant $1$, in $\R^3$.

\begin{prop}
Let $V:T\widetilde{S_g}\to T\mathbb{X}$ be as in the previous section and let $x\in \widetilde{S_g}$. The plane $V(T_x\widetilde{S_g})\subset \mathbb{X}$ is tangent to a reducible $\mathbb{H}^2$.
\end{prop}

\begin{proof}
	We fix a unitary basis of the Higgs bundle $\mathcal{E}_x$ over $x$ for the harmonic metric, adapted with the line decomposition $\mathcal{E}_x=\mathcal{K}^{1/2}_x\oplus \mathcal{O}_x\oplus \mathcal{K}^{-1/2}_x$. In this unitary basis the plane $V(T_x\widetilde{S_g})$ can be identified with  the plane of Hermitian matrices :
$$\begin{pmatrix}
0 & 0 & \overline{\tau}\\
0 & 0 & 0\\
\tau & 0 & 0
\end{pmatrix},\,\tau \in \mathbb{C}$$

In such a basis the antilinear involution whose fixed point determines the real locus is given by $(z_1,z_2,z_3)\mapsto (\overline{z_3},\overline{z_2},\overline{z_1})$.
	So in order to obtain the corresponding real symmetric matrices, we apply the change of variable $z_1=\frac{1}{\sqrt{2}}(x_1+ix_3)$, $z_2=x_2$, $z_3=\frac{1}{\sqrt{2}}(x_1-ix_3)$. In this new basis, the plane $V(T_x\widetilde{S_g})$ is the plane of symmetric real matrices : 

$$\begin{pmatrix}
\sqrt{2}\text{Re}(\tau) & 0 & \sqrt{2}\text{Im}(\tau)\\
0 & 0 & 0\\
\sqrt{2}\text{Im}(\tau) & 0 & -\sqrt{2}\text{Re}(\tau)
\end{pmatrix},\,\tau \in \mathbb{C}$$

This plane is the tangent plane to $\mathcal{H}_0$ from above, as claimed.

\end{proof}

We now describe a projection $p:\mathcal{F}_{1,2}\to \mathcal{H}$. We will however mostly need to use the more concrete description of this map, detailed in the next subsection.

\medskip

Let $\mathcal{H}\subset \mathbb{X}$ be a reducible hyperbolic plane. The nearest point projection defines a map $p_0:\mathbb{X}\to \mathcal{H}$. This map extends in some sense to a map $p:\mathcal{F}_{1,2}\to\overline{\mathcal{H}}= \mathcal{H}\cup \partial \mathcal{H}$, where $\partial \mathcal{H} \simeq S^1$ is the visual boundary of $\mathcal{H}$. One can define this extension in terms of minimum of Busemann functions as in \cite[{Theorem 7.3}]{Dav23}. The fact that this projection extends well is due to the fact that the tangent vectors to $\mathcal{H}$ are regular (their Cartan projection do not lie in any wall of the model Weyl chamber).

\medskip

Let us fix a basis of $\mathbb{R}^3$ and a dual basis of $\left(\mathbb{R}^3\right)^*$. A flag $f\in \mathcal{F}_{1,2}$ is parametrized by a pair $[x]=[x_1:x_2:x_3]\in \mathbb{P}\left(\R^3\right)$ and $[y]=[y_1:y_2:y_3]\in \mathbb{P}\left( \left(\R^3\right)^*\right)$ such that $x\cdot y=x_1y_1+x_2y_2+x_3y_3=0$.
The symmetric space $\mathbb{X}$ associated with $\SL(3,\R)$ can be identified with the space of symmetric positive matrices with determinant one. 

We identify $\mathcal{F}_{1,2}$ with the $\SL(3,\R)$-orbit in the visual boundary of $\mathbb{X}$ in such a way that the Busemann function associated with $f=([x],[y])$ is equal, for some scaling of the metric on $\mathbb{X}$, for $O \in\mathbb{X}$ to:

\begin{equation}\label{eq:BusemannFunction}
b_{f,O}:X\in \mathbb{X}\mapsto \log\left(\frac{x^tXx}{x^tOx}\right)+\log\left(\frac{y^tX^{-1}y}{y^tO^{-1}y}\right).
\end{equation}

Note that the $\SL(3,\R)$-orbit in the visual boundary that we identified with $\mathcal{F}_{1,2}$ is coincidentally also the orbit containing $\partial \mathcal{H}$.
 
The vector field $\mathrm{v}_f=-\nabla b_{f,O}$ is what we will call the vector field of vectors \emph{pointing towards $f\in \mathcal{F}_{1,2}$}.  

\begin{prop}[{\cite[{Lemma 7.2 and Theorem 7.11}]{Dav23}}]
The function $b_{f,O}$ is convex and Lipschitz, and $e^{\lambda b_{f,O}}$ is strictly convex on $\mathcal{H}$ for some $\lambda>0$. 

In particular it admits either a unique critical point that is a global minimum $x\in \mathbb{H}$, or there is a unique point $a\in \partial \mathcal{H}$ in the visual boundary of $\mathcal{H}$ such that $b_{f,O}$ does not go to $+\infty$ on some and hence any geodesic rays converging to $a$.

In the second case $b_{f,O}$ actually goes to $-\infty$ on the geodesic rays converging to $a$.  
\end{prop}

Let $p:\mathcal{F}_{1,2}\to \mathcal{H}\cup \partial \mathcal{H}$ be the map that associates to $f\in \mathcal{F}_{1,2}$ either the minimum $x\in \mathcal{H}$ of $b_{f,O}$ or the unique point $a\in \partial \mathcal{H}$ such that $b_{f,O}$ converges to $-\infty$ for some and hence all geodesic rays converging to $a$. It is well defined because of the previous proposition.

\subsection{Description of the fibration}

We now give a description of the preimage by $p$ of a point and a geodesic. We show that the pre-images of half-spaces are multicones with $\mathcal{C}^1$ boundary.

\medskip

Recall that $\partial \mathcal{H}$ can be identified with a circle inside $\mathcal{F}_\Delta\simeq \mathcal{F}_{1,2}$. Recall that the thickening of a flag was introduced in Section \ref{subsec:Multicone}.

\begin{prop}
\label{prop:Image boundary fibration}
Let $a\in \partial \mathcal{H}$. The pre-image $p^{-1}\left(\lbrace a\rbrace\right)$ is equal to the thickening of $K_a\subset \mathcal{F}_{1,2}$, where $a$ is identified with an element of $\mathcal{F}_\Delta$.
\end{prop}

Note that $K_a\cap K_b=\emptyset$ for $a\neq b\in \partialvis \mathcal{H}$ since any two such flags are transverse.

\begin{proof}
The asymptotic slope of a Busemann function \eqref{eq:BusemannFunction} along a geodesic ray is given by the Tits angle \cite[{Lemma 4.8}]{Dav23}. The set $K_a$ is exactly the set of elements of $\mathcal{F}_{1,2}$ whose Tits angle with the end point of any geodesic ray converging to $a\in \partial \mathcal{H}$ is less than $\frac{\pi}{2}$, which is therefore also the set of flags whose associated Busemann function does not diverge to $+\infty$ on rays converging to $a$. In other words, these are exactly the elements of $p^{-1}\left(\lbrace a\rbrace\right)$. 
\end{proof}

We can now prove the following:

\begin{prop}
The map $p$ is continuous.
\end{prop}

\begin{proof}
This map is continuous on the pre-image of $\mathcal{H}$ by \cite[{Theorem 7.3}]{Dav23}. By Proposition \ref{prop:Image boundary fibration}, since the compact set $K_a$ varies continuously in $a$, $p$ is continuous on $\partial \mathcal{H}$. If we now consider a sequence $(f_n)$ in $\mathcal{F}_{1,2}$ converging to $f$ such that $p(f_n)\in \mathcal{H}$ converges to $a\in \partial \mathcal{H}$, we know that the Buseman function $b_{f,O}$ for $O\in \mathcal{H}$ is bounded from above on the geodesic ray between $O\in \mathcal{H}$ and $a\in \mathcal{F}_{1,2}$, therefore $p(f)=a$.
\end{proof}

We say that a \emph{conic in $\mathcal{F}_{1,2}$} is the set of pairs $f=(\ell, H)$ for some fixed non-degenerate conic $C$ in $\mathbb{RP}^2$ where $\ell$ lies on the conic and the projective line defined by $H$ is tangent to the conic.  

\begin{prop}
\label{prop:Conic Fiber}
Let $X\in  \mathcal{H}$. The pre-image $p^{-1}\left(\lbrace X\rbrace\right)$ is equal to a conic in $\mathcal{F}_{1,2}$ corresponding to a conic $C(X,\mathcal{H})$ in $\mathbb{RP}^2$. More precisely if $\mathcal{H}={\mathcal{H}_0}$ and $X=\text{Id}$, then in coordinates this set can be written as:
$$\lbrace ([\cos(\theta):1:\sin(\theta)]  ,[-\cos(\theta):1:-\sin(\theta)] )\mid \theta\in [0, 2\pi]\rbrace.$$

If $\mathcal{H}=\mathcal{H}_0$ and $X=\Id$, this conic $C(\Id, \mathcal{H}_0)$ is defined in coordinates by $x_1^2-x_2^2+x_3^2=0$.
\end{prop}

\begin{proof}
Let us assume that $\mathcal{H}=\mathcal{H}_0$ and $X=\Id$ up to acting by $\SL(3,\R)$. Given a flag $f$, one has $p(f)=\Id$ if and only if $\Id$ is a critical point for $b_{f,O}$ \eqref{eq:BusemannFunction}. If $f$ in coordinates is equal to $([x_1:x_2:x_3],[y_1:y_2:y_3])$, this can be written as:
$$\frac{y_1^2-y_3^2}{y_1^2+y_2^2+y_3^2}=\frac{x_1^2-x_3^2}{x_1^2+x_2^2+x_3^2},$$
$$\frac{2y_1y_3}{y_1^2+y_2^2+y_3^2}=\frac{2x_1x_3}{x_1^2+x_2^2+x_3^2}.$$

We can assume that $(x_1^2+x_2^2+x_3^2)=(y_1^2+y_2^2+y_3^2)=1$ up to scaling by a positive element. Therefore:
$$y_1^2-y_3^2=x_1^2-x_3^2  $$
$$y_1y_3=x_1x_3 $$

In particular the squares of the complex numbers $x_1+ix_3$ and $y_1+iy_3$ are equal so up to multiplying $y$ by $-1$ one can assume $x_1=-y_1$ and $x_3=-y_3$. Note that moreover this implies that $x_2^2=y_2^2$ since the sum of the squares are equal. The fact that $x\cdot y=0$ implies that $x_1^2+x_3^2= x_2y_2$. In particular $x_2=y_2$ so $f=([x_2\cos(\theta):x_2:x_2\sin(\theta)]  ,[-x_2\cos(\theta):x_2:-x_2\sin(\theta)] )$ for some $\theta\in[0, 2\pi]$. This implies the desired result. The corresponding conic in $\mathbb{RP}^2$ is the one associated with the equation $x_1^2-x_2^2+x_3^2=0$. 

%
%
\end{proof}

Given a non-zero vector $v\in T\mathcal{H}$, we consider the orthogonal geodesic $\gamma\subset \overline{\mathcal{H}}$. The complement $\overline{\mathcal{H}}\setminus \gamma$ is the union of two open half spaces. Let $H_v\subset \overline{\mathcal{H}}$ be the open half-space containing $v$. We define $ U^\mathcal{H}_v=p^{-1}(H_v)$. 

\begin{lem}
The set $U^\mathcal{H}_v$ is a multicone. Moreover $\partial U^\mathcal{H}_v$ is a $\mathcal{C}^1$ hypersurface and $\partial U^\mathcal{H}_v\cap p^{-1}(\mathcal{H})$ is dense in $\partial U^\mathcal{H}_v$.
\end{lem}
%

\begin{proof}

If we denote by $f^+, f^-$ the full flags corresponding respectively to the endpoint of the geodesics with initial derivatives $v$ and $-v$, one has $K_{f^+}\subset U^\mathcal{H}_v$ and $K_{f^-}\cap\overline{U^\mathcal{H}_v}=\emptyset$. Hence $U^\mathcal{H}_v$ is a multicone.

\medskip

Assume $\mathcal{H}=\mathcal{H}_0$ and let $\gamma$ be the geodesic such that in coordinates the endpoints $\gamma^+,\gamma^-\in \partial \mathcal{H}$ of $\gamma$  are the flags:
$$\gamma^+=([1:0:0],[0:0:1]),$$
$$\gamma^-=([0:0:1],[1:0:0]).$$

\medskip

The set $\partial U^\mathcal{H}_v\cap p^{-1}(\mathcal{H})$ admits the following smooth parametrization:
$$\lbrace ([\lambda\cos(\theta):1:\lambda^{-1}\sin(\theta)]  ,[-\lambda^{-1}\cos(\theta):1:-\lambda\sin(\theta)] )\mid \theta\in [0, 2\pi], \lambda>0\rbrace.$$

The boundary $\partial U_v^\mathcal{H}$ decomposes as:
$$p^{-1}(\partial H_v)=p^{-1}(\gamma\cap\mathcal{H})\cup p^{-1}(\lbrace \gamma^+\rbrace)\cup p^{-1}(\lbrace \gamma^-\rbrace).$$
It remains to show that $\partial U^\mathcal{H}_v$ is smooth in a neighborhood of $K_{\gamma^+}$ and $K_{\gamma^-}$.

\medskip

The thickening $K_{\gamma^+}$ decomposes as a wedge or circles:
$$K_{\gamma^+}=S_1\cup S_2\cup \lbrace \gamma^+\rbrace,$$
$$ S_1=\lbrace ([a:1:0],[0:0:1]) \mid a\in \mathbb{R}\rbrace,$$
$$S_2=\lbrace ([1:0:0],[0:1:a]) \mid a\in \mathbb{R}\rbrace, $$
$$\gamma^+= ([1:0:0],[0:0:1]).$$ 

 We first study a neighborhood of $S_1$, the case of $S_2$ being completely analog. Consider the chart of the space of flags with coordinates $x_1=1$ and $y_2=1$. In this chart the elements of $\partial U^\mathcal{H}_v\cap p^{-1}(\mathcal{H})$ can be written for $\lambda>0$ and $\cos(\theta)\neq 0$:
$$([1:\lambda^{-1}\cos^{-1}(\theta):\lambda^{-2}\tan(\theta)]  ,[-\lambda^{-1}\cos(\theta):1:-\lambda\sin(\theta)] ).$$

The union of this set and $S_2$ is parametrized in coordinates as follows::
$$\left\lbrace\left(\left[1:x_2:-\frac{x_2^3y_3}{1+x_2^2y_3^2}\right]  ,\left[-\frac{x_2}{1+x_2^2y_3^2}:1:y_3\right] \right) \right\rbrace.$$

Indeed one can take freely any non-zero value for $x_2$ and any value for $y_3$ by choosing $\theta$ and $\lambda$, and the other coordinates are determined. When $x_2=0$ this parametrizes $S_2$.

\medskip

This is the graph of a smooth map, hence it is a smooth submanifold.

\medskip

 We now study a neighborhood of $\gamma^+$. Consider the chart of the space of flags with coordinates $x_1=1$ and $y_3=1$. In this chart the elements of $\partial U^\mathcal{H}_v\cap p^{-1}(\mathcal{H})$ can be written for $\lambda>0$ and $\cos(\theta)\neq 0$:
$$([1:\lambda^{-1}\cos^{-1}(\theta):\lambda^{-2}\tan(\theta)]  ,[-\lambda^{-2}\cot(\theta):\lambda^{-1}\sin^{-1}(\theta):1] ).$$

The union of this set with $S_1$ and $S_2$ minus two points is parametrized similarly in coordinates as follows:
$$\left\lbrace\left(\left[1:x_2:\frac{x_2^3y_2}{x_2^2+y_2^2}\right]  ,\left[-\frac{y_2^3x_2}{x_2^2+y_2^2}:y_2:1\right] \right) \right\rbrace.$$

Indeed one can take freely any non-zero value for $x_2$ and $y_2$ by choosing $\theta$ and $\lambda$, and the other coordinates are determined. When $x_2=0$ this parametrizes $S_2\cup\gamma^+$ minus a point and when $y_2=0$ this parametrizes $S_1\cup\gamma^+$ minus a point. Note that this expression naturally extends for $x_2=y_2=0$ by continuity.

\medskip

This is the graph of a function that is $\mathcal{C}^1$ (but not $\mathcal{C}^2$), so we have a $\mathcal{C}^1$ submanifold.

\medskip

Finally note that every point of $K_{\gamma^+}$ and $K_{\gamma^-}$ can be approximated by points in $\partial U^\mathcal{H}_v\cap p^{-1}(\mathcal{H})$. The only two points for which this surface is not $\mathcal{C}^2$ are $\gamma^+$ and $\gamma^-$.
\end{proof}

\begin{rem}
The surface $\partial U^\mathcal{H}_v$ is a genus two orientable surface. Indeed it is co-oriented in the orientable manifold $\mathcal{F}_{1,2}$, and can be obtained as the gluing of a cylinder along the disjoint union of two wedges of circles.
\end{rem}

\section{Geometric structure}
\label{sec:Geom}

In this section with use the parallel distribution of planes in the symmetric space $\mathbb{X}=\SL(3,\R)/SO(3,\R)$ and the construction from Section \ref{sec:MulticonesConstruction} to construct a family of multicones associated to any representation in $\mathcal{S}_X^\R$. We use them to prove the Anosov property and to describe an associated geometric structure. 

\subsection{Nestedness of the multicones.}
\label{subsec:Nestedness}
We perform the technical computations needed for the proof of Theorem \ref{thm:main}: we use the bound given by the maximum principle to show that the multicones we defined are nested, in order to apply the criterion of Theorem \ref{thm:CriterionAnosov}.

\medskip

Recall that we consider a representation $\rho:\pi_1(S_g)\to \SL(3,\R)$ associated to a Higgs bundle in $\mathcal{S}_X^\R$, and we constructed an equivariant map $V:T\widetilde{S_g}\to T\mathbb{X}$ that defines a parallel distribution of planes in the symmetric space.

\medskip

Given a non-zero vector $v\in T_x\widetilde{S_g}$ we associate the multicone $U_v=U^{\mathcal{H}_x}_{V(v)}$, where $\mathcal{H}_x$ is the reducible hyperbolic plane spanned by $V(T_x\widetilde{S_g})$.

\begin{thm}
\label{thm:Nestedness}
Let $\Phi$ be the geodesic flow on $\widetilde{S_g}$ for the metric $h_1^{-2}$. For all $t>0$ and all non-zero vectors $v\in T\widetilde{S_g}$, $\overline{U_{\Phi_t(v)}}\subset U_v$. More precisely the $\mathcal{C}^1$ hypersurfaces $(\partial U_{\Phi_t(v)})_{t\in \R}$ define locally a fibration. 
\end{thm}

\begin{proof}

Let $x\in \widetilde{S_g}$ and $v\in T_x\widetilde{S_g}$. We first use the parallelism property to understand how the multicone vary at first order.

The first order variation of $V(v)$ and $V(T_x\widetilde{S_g})$ when $(x,v)$ moves along the geodesic flow $\Phi$ for the metric $h_1^{-2}$ is given by the infinitesimal transvection associated to $H=\mathrm{d}h(v)\in T\mathbb{X}$ because of Proposition \ref{prop:ParallelPlanes}. Hence the first order variation of the hypersurface $\partial U_v$, which depends smoothly and $G$-equivariantly on $V(v)$ and $V(T_x\widetilde{S_g})$ is also given by the infinitesimal transvection associated to $H$ when $v$ moves along the geodesic flow $\Phi$.

\medskip

 The $\mathcal{C}^1$ hypersurface $\partial U_v$ is co-oriented, hence its normal bundle $T\mathcal{F}_{1,2}/T\partial U_v$ admits some non-vanishing linear form $\eta$ positive on vectors pointing inside $U_v$.  The infinitesimal transvection associated to $H$ defines a vector field $W_v$ on $\mathcal{F}_{1,2}$ defined at $f\in \mathcal{F}_{1,2}$ as the derivative at $t=0$ of $t\mapsto e^{tH}\cdot f$.

\medskip

If $\eta(W_v)>0$ on $\partial U_v$ for all $v\in T\widetilde{S_g}$, the hypersurfaces $\partial U_{\Phi_t(v)}$ are locally nested when $t$ varies. As a consequence $\overline{U_{\Phi_t(v)}}\subset U_v$ for all $t>0$, and the proof would be done.

\medskip

We will only define a one form of the normal bundle $\eta$ on a dense open subset of $\partial U_v$ (the whole surface minus two wedges of two circles). We will construct this one form so that for any Riemannian metric on $\mathcal{F}_{1,2}$ the norm of $\eta$ is bounded on this dense set. We then only need to show that for some $\epsilon>0$ one has on this dense subset:
\begin{equation}\label{eq:positivity}
\eta(W_v)>\epsilon
\end{equation}

To conclude the proof we now construct the bounded one-form $\eta$ positive on vectors pointing towards $U_v$ and show the condition \eqref{eq:positivity}. 

\medskip

For that purpose we fix $x\in \widetilde{S_g}$ and $v\in T_x\widetilde{S_g}$. We fix a unitary basis of the Higgs bundle $\mathcal{E}_x$ over $x$ for the harmonic metric , adapted with the line decomposition $\mathcal{E}_x=\mathcal{K}^{1/2}_x\oplus \mathcal{O}_x\oplus \mathcal{K}^{-1/2}_x$. The real form associated with the harmonic metric is equal in this coordinates to:
 $$(z_1,z_2,z_3)\mapsto (\overline{z_3},\overline{z_2},\overline{z_1}).$$
 
\begin{rem}

From now on we work with matrices with complex entries, but which correspond to real endomorphisms for this non-standard real form. Hence what we call symmetric matrices from now on are hermitian matrices that are real endomorphisms for this modified real form. 
 
\end{rem}

In such coordinates, up to acting by $(z_1,z_2,z_3)\mapsto (wz_1,z_2,\overline{w}z_3)$ for $w$ a unit complex number, the symmetric matrix associated with $\mathrm{d}h(v)$ can be written as:
 $$ H=\begin{pmatrix}
0 & \beta & 1\\
\overline{\beta} & 0 & \beta\\
1 & \overline{\beta} & 0
\end{pmatrix}\simeq \mathrm{d}h(v).$$

Here $\beta$ is a complex number whose modulus is equal to $\Vert t\Vert/\Vert \tau \Vert$. Note that the maximum principle from Proposition \ref{prop:Maximal Principle} implies that $|\beta|<1$.

%

\medskip

We now proceed in constructing this one form $\eta$ and checking the positivity assumption. In this coordinates, the plane $V(T_x\widetilde{S_g})$ tangent to the reducible hyperbolic plane $\mathcal{H}_x$ is generated by the symmetric matrices:

 $$ H_0=\begin{pmatrix}
0 & 0 & 1\\
0 & 0 & 0\\
1 & 0 & 0
\end{pmatrix}\simeq V(v),$$

 $$ H_0^\perp=\begin{pmatrix}
0 & 0 & -i\\
0 & 0 & 0\\
i & 0 & 0
\end{pmatrix}\simeq V(J(v)).$$

\medskip

\begin{figure}
\begin{center}
\begin{tikzpicture}[scale=1.5]
\tdplotsetmaincoords{70}{90}
\tdplotsetrotatedcoords{0}{0}{0};
\draw[ thick,
    tdplot_rotated_coords] (0,2,0) arc (90:270:2);

\tdplotsetmaincoords{0}{0}
\fill (0,0,0) circle (.025);
\node[above left] () at (0,0,0) {$x$};
\node[above left] () at (0,0,-1) {$x_d$};
\node[left] () at (-2,0,0) {$\mathcal{H}_x$};
\node[above right] () at (0,0,-1.7) {$\gamma$};
\draw[red,thick] (0,0,1.7) -- (0,0,-1.7);
\fill (0,0,-1) circle (.025);
\draw[blue, very thick, ->] (0,0,0) -- (0,0,.5);
\draw[blue,very thick, ->] (0,0,0) -- (.5,0,0);
\draw[very thick, ->] (0,0,0) -- (.5,.3,0);
\node[right] () at (.5,0,0) {$H_0$};
\node[below] () at (0,0,.5) {$H_0^\perp$};
\node[right] () at (.5,.3,0) {$H$};

\draw[thick,
    tdplot_rotated_coords] (0,2,0) arc (90:-90:2);

\end{tikzpicture}
\end{center}
\caption{Illustration of the reducible hyperbolic plane $\mathcal{H}_x$ in the symmetric space}
\label{fig:FibrationDoD}
\end{figure}
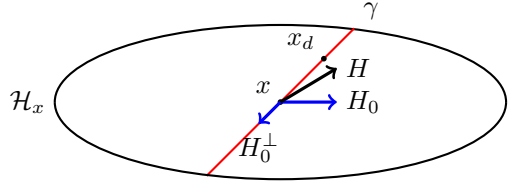

Let $f=(f^1,f^2)\in \mathcal{F}_{1,2}$ be a full flag in $\C^3$, coming from a real flag. To such a flag we associate $\pi=(\pi_{i,j})_{i,j}$ a rank one matrix with kernel $f^2$ and image $f^1$. Recall that $p:\mathcal{F}_{1,2}\to \overline{\mathcal{H}_x}$ is the extension of the nearest point projection. 

\medskip

We now check that $p^{-1}\left(\lbrace h(x)\rbrace \right)$ is exactly the circle that can be parametrized as the projectors $\pi$ having the following form, for all $z\in \C$ with $|z|=1$:
\begin{equation}
\label{eq:defPi}
\pi=\frac{1}{4}\begin{pmatrix}
-1 & \sqrt{2} z & -z^{2} \\
-\sqrt{2} \overline{z} & 2 & -\sqrt{2} z \\
-\overline{z}^{2} & \sqrt{2} \overline{z} & -1
\end{pmatrix}.
\end{equation}

The corresponding flag in coordinates is given by $([z: \sqrt{2}: \overline{z}], [-z: \sqrt{2}: -\overline{z}])\in \mathbb{P}(\R^3)\times \mathbb{P}\left(\left(\R^3\right)^*\right)$. These coordinates describe a flag since $-z\overline{z}+2-z\overline{z}=0$, that is real for the non-standard real form. 

\medskip

The matrix $\mathrm{v}_f=\pi\pi^*-\pi^*\pi$ corresponding to the the vector in $T_{h(x)}\mathbb{X}$ at our basepoint $h(x)$ pointing towards the flag $f$ is equal to $\pi\pi^*-\pi^*\pi$, which is equal to:
$$\begin{pmatrix}
0 & \frac{\sqrt{2}}{2} z & 0 \\
\frac{\sqrt{2}}{2} \overline{z}& 0 & \frac{\sqrt{2}}{2} z \\
0 & \frac{\sqrt{2}}{2} \overline{z} & 0
\end{pmatrix}.$$

This matrix is collinear with the gradient of the Busemann function associated with $f$, and it is orthogonal to the plane spanned by $H_0$ and $H_0^\perp$. Hence we see that the circle of flags obtained for all unit complex $z$ belongs to $p^{-1}\left(\lbrace x\rbrace \right)$. We already know that $p^{-1}\left(\lbrace x\rbrace \right)$ is a circle, so as we claimed these flags form exactly the fiber $p^{-1}\left(\lbrace h(x)\rbrace \right)$.

\medskip

 We have an identification $f\mapsto \mathrm{v}_f\in T_{h(x)}\mathbb{X}$ of the real flag manifold $\mathcal{F}_{1,2}$ with a closed set of Hermitian matrices. The complex linear form that associates to a Hermitian matrix $m=(m_{i,j})_{i,j}$ the coefficient $m_{1,3}=\overline{m_{3,1}}$ defines a complex-valued one form $\alpha$ on $\mathcal{F}_{1,2}$ whose kernel is tangent to the circle $p^{-1}\left(\lbrace h(x)\rbrace \right)\subset \mathcal{F}_{1,2}$.

\medskip

Let $\gamma\subset \mathcal{H}$ be the geodesic through $h(x)$ orthogonal to $V(v)$. The tangent vector of $\gamma$ at $h(x)$ corresponds to the Hermitian matrix $H_0^\perp$.

\medskip

Note that $\partial U_v$ decomposes as $p^{-1}(\gamma)\cup K_{f^+}\cup K_{f^-}$ where $f^+,f^-$ are the flags associated with the endpoints of $\gamma$. The surface $p^{-1}(\gamma)$ is a cylinder that is the union of the circles $e^{dH_0^\perp}\cdot p^{-1}\left(\lbrace h(x)\rbrace \right)$ for $d\in \mathbb{R}$. 

\medskip

Let $\alpha_0$ be function on $\mathcal{F}_{1,2}$ that associates to a flag $f$ the derivative at $t=0$ of $t\mapsto \alpha(e^{tH_0^\perp}\cdot f)$. The tangent space of the surface $p^{-1}(\gamma)$ is equal on $p^{-1}(\lbrace h(x)\rbrace)$ to the set of vectors $w\in T\mathcal{F}_{1,2}$ such that $\alpha(w)$ is collinear with $\alpha_0$. In particular the imaginary part $\eta_0=\text{Im}(\overline{\alpha_0}\alpha)$ defines a real one form on the normal bundle of the surface $p^{-1}(\gamma)$ in $\mathcal{F}_{1,2}$.

\medskip

 Let $x_d=e^{dH_0^\perp}\cdot h(x)$ for $d\in \mathbb{R}$. We define $\eta$ to be equal to $\frac{e^{-dH_0^\perp}\cdot \eta_0}{\cosh(d)^2}$ on $p^{-1}(x_d)$. This form is defined on $p^{-1}(\gamma)$ which is dense in $\partial U_v$.
 
 \begin{lem}
 \label{lem:small Lemma bounded}
 For any Riemannian metric on $\mathcal{F}_{1,2}$, the norm of $\eta$ on $p^{-1}(\gamma)$ is bounded. 
 \end{lem}
 
 \begin{proof}[{Proof of Lemma \ref{lem:small Lemma bounded}}]
 The differential of the action of $e^{dH_0^\perp}$ on the flag manifold $\mathcal{F}_{1,2}$ has norm at most $O(e^{2|d|})$. Hence the one-form $\frac{e^{dH_0^\perp}\cdot \eta_0}{\cosh(d)^2}$ has bounded norm.
 \end{proof}

We now do the main computation:

 \begin{lem}
 \label{lem:Main Computations}
Let $\epsilon=\frac{1-|\beta|}{2}>0$ and let $W_v$ as before be the vector field on $\mathcal{F}_{1,2}$ that is equal at $f\in \mathcal{F}_{1,2}$ to the derivative at $t=0$ of $t\mapsto e^{tH}\cdot f$. On $p^{-1}(\gamma)$ one has:
$$\eta(W_v)\geq\epsilon.$$
 \end{lem}
 
 Once this lemma is proven, the last thing that we need to check is that $\eta$ is positive on the correct side of $\partial U_v$, i.e. on vectors pointing towards $U_v$. but this is equivalent to having \eqref{eq:positivity} in the case $\beta=0$, so it also follows from Lemma \ref{lem:Main Computations}.
 
 \end{proof}
  
We now prove the last Lemma that we used to conlude the proof of Theorem \ref{thm:Nestedness}.
 
 \begin{proof} [{Proof of Lemma \ref{lem:Main Computations}}]
Let $e^{tH_0^\perp}\cdot f\in p^{-1}(x_d)$ for some $d\in \R$ and $f\in p^{-1}\left(\lbrace h(x)\rbrace \right)$. Let $\pi$ be the rank one matrix associated with $f$ as in \eqref{eq:defPi}, for some complex number $z$ such that $|z|^2=1$. The action of $\SL(3,\R)$ on the space of rank one projectors is by conjugation, and the derivative of the adjoint action of $e^{tH}$ on $\pi$ is equal to $[H, e^{dH^\perp_0}\cdot \pi]=e^{dH^\perp_0}\cdot[H',\pi]$ where $H'=e^{-dH_0^\perp}He^{dH_0^\perp}$.

\medskip

The derivative of the norm one hermitian matrix $\mathrm{v}_f$ pointing towards $f$ under the action of $(e^{tH})_{t\in \mathbb{R}}$ is equal to $e^{dH^\perp_0}\cdot\left([H',\pi]\pi^*+\pi[H',\pi]^*-[H',\pi]^*\pi-\pi^*[H',\pi]\right)$.

\medskip

Hence we need to show that for some $\epsilon>0$ depending only on $|\beta|$:
$$\text{Im}\left(\alpha(M_1)\overline{\alpha(M_2)}\right)\geq\epsilon \cosh(d)^2 .$$
In this expression:
$$M_1=[H',\pi]\pi^*+\pi[H',\pi]^*-[H',\pi]^*\pi-\pi^*[H',\pi],$$
$$M_2=[H_0^\perp,\pi]\pi^*+\pi[H_0^\perp,\pi]^*-[H_0^\perp,\pi]^*\pi-\pi^*[H_0^\perp,\pi].$$

To compute this expression note that:

$$e^{dH_0^\perp}=\begin{pmatrix}
\cosh(d) & 0 & i\sinh(d)\\
0 & 1 & 0\\
-i\sinh(d) & 0 & \cosh(d)
\end{pmatrix}.$$

And therefore :

\newcommand*{\Scale}[2][4]{\scalebox{#1}{$#2$}}%

$$H'=\Scale[1]{\begin{pmatrix}
2 i \cosh\left(d\right) \sinh\left(d\right) & \beta \cosh\left(d\right) + i  \overline{\beta} \sinh\left(d\right) & \cosh\left(d\right)^{2} + \sinh\left(d\right)^{2} \\
\overline{\beta} \cosh\left(d\right) + i \beta \sinh\left(d\right) & 0 & \beta \cosh\left(d\right) - i  \overline{\beta} \sinh\left(d\right) \\
\cosh\left(d\right)^{2} + \sinh\left(d\right)^{2} & \overline{\beta} \cosh\left(d\right) - i \beta \sinh\left(d\right) & -2 i \cosh\left(d\right) \sinh\left(d\right)
\end{pmatrix}}.$$
One can directly compute the $(1,3)$-coefficient of the matrices $M_1$ and $M_2$. We obtain:
$$\alpha(M_1)=\left(3-\overline{z}^4\right)\frac{\cosh(d)^2+\sinh(d)^2}{4}-\sqrt{2}i\beta\overline{z}\sinh(d),$$
$$\alpha(M_2)=\left(3+\overline{z}^4\right)\frac{i}{4}$$

As a consequence:
$$\text{Im}\left(\alpha(M_1)\overline{\alpha(M_2)}\right)\geq \frac{\cosh(d)^2+\sinh(d)^2}{2}-\sqrt{2}|\beta|\sinh(|d|) .$$

Note that $\cosh(d)^2+\sinh(d)^2=1+2\sinh(|d|)^2\geq 2\sqrt{2}\sinh(|d|)$ therefore:
$$\text{Im}\left(\alpha(M_1)\overline{\alpha(M_2)}\right)\geq \frac{\cosh(d)^2+\sinh(d)^2}{2}(1-|\beta|).$$

We see that $\epsilon=\frac{1-|\beta|}{2}$ works.
 \end{proof}
 

We can now prove our main Theorem.

\begin{thm}
\label{thm:mainText}
Representations corresponding to the Higgs bundles $\mathcal{S}_X^\R$ are Anosov.
\end{thm}

\begin{proof}

We consider the metric $h_1^{-2}$ induced by $H$ on $X$. We apply the criterion of Theorem  \ref{thm:CriterionAnosov} to the familly of multicones $(U_v)$, which satisfy the nestedness property required because of Theorem \ref{thm:Nestedness}, for the geodesic flow associated with the metric $h_1^{-2}$. Hence the underlying representation $\rho$ is Anosov.
\end{proof}

\subsection{Fibration of a domain of discontinuity}

Given a representation $\rho$ associated with a Higgs bundle in $\mathcal{S}_X^\R$, the parallel distribution of planes describes a fibration of a domain of disconitnuity for $\rho$ in $\mathcal{F}_{1,2}$.

\medskip

Let $\Omega_\rho=\mathcal{F}_{1,2}\setminus \bigcup_{\zeta\in \partial \Gamma_g}K_{\xi_\rho(\zeta)}$. This domain is a cocompact domain of discontinuity for $\rho$, beauce $\rho$ is Anosov. This domain is a particular examples of domains constructed by Guichard-Wienhard \cite{GWDoD} and Kapovich-Leeb-Porti \cite{KLP}. Barbot showed that the quotient of this domain by the representation is a circle bundle over the surface $S_g$ for any non-Hitchin Anosov representation $\rho:\Gamma_g\to \SL(3,\R)$ \cite{Bar10}. 

\medskip

A reducible copy $\mathcal{H}$ of the hyperbolic plane in $\mathbb{X}$ together with a point $o\in\mathcal{H}$ define a conic $C(x,\mathcal{H})$ by Proposition \ref{prop:Conic Fiber}.

\begin{thm}
\label{thm:Fibration}
The conic of flags associated to the conics $(C(h(x), \mathcal{H}_x))_{x\in \widetilde{S_g}}$ form a fibration of the domain of discontinuity $\Omega_\rho$.
\end{thm}

\begin{proof}
Let $v\in T_x\widetilde{S_g}$. Theorem \ref{thm:Nestedness} implies that the hypersurfaces $(\partial U_{\Phi_t(v)})_{t\in \R}$ define locally a fibration. However if $x_t\in \widetilde{S_g}$ is the basepoint of $\Phi_t(x)$, then $C(h(x), \mathcal{H}_x)\subset \partial U_{\Phi_t(v)}$, therefore the conics $(C(h(x_t), \mathcal{H}_{x_t}))_{t\in \R}$ define a fibration of a hypersurface in $\mathcal{F}_{1,2}$. This holds for all tangent vectors $v$, hence the projection on the the second factor of $\lbrace (x,f)\mid x\in \widetilde{S_g},\,f\in C(h(x), \mathcal{H}_x)\rbrace$ is an immersion, which is proper so it locally defines a fibration.

\medskip

We now check that for $x\neq y$, the conics $C(h(x), \mathcal{H}_x)$ and $C(h(y), \mathcal{H}_y)$ are disjoint. By taking a geodesic from $x$ to $y$ one can find $v\in T_x\widetilde{S_g}$ such that $\Phi_t(v)$ is based at $y$. One has $\overline{U_{\Phi_t(v)}}\subset U_v$ and hence $\partial U_{\Phi_t(v)}\cap \partial U_{v}=\emptyset$ and in particular $C(h(x), \mathcal{H}_x)\cap C(h(y), \mathcal{H}_y)=\emptyset$.

\medskip

Now we have a domain in $\mathcal{F}_{1,2}$ that is fibered by conics. We check that this domain is $\Omega_\rho$. Given any $\zeta\in\partial \Gamma$ and any $x\in \widetilde{S_g}$, we find a geodesic for the metric $h_1^{-2}$ from $x$ to $\zeta$. By the last part of Theorem \ref{thm:CriterionAnosov}, $K_{\xi_\rho(\zeta)}\subset U_v$, and therefore $C(h(x), \mathcal{H}_x)$ is disjoint from $K_{\xi_\rho(\zeta)}$. Note that the metric associated to $h_1^{-2}$ is negatively curved, see Remark \ref{rem:NegativeCurved}. 

Hence the domain fibered by these conics is contained in $\Omega_\rho$. The quotient of this domain by $\rho(\Gamma_g)$ is a circle bundle over $\widetilde{S_g}$, hence it is compact and therefore closed in $\Omega_\rho/\rho(\Gamma_g)$. Is is also open, so it is a connected component. 

The domain $\Omega_\rho$ is connected because $\rho$ can be deformed by a path of Anosov representation into a reducible representation, for which the domain of discontinuity is connected. Hence it is exactly $\Omega_\rho$ that is fibered by the conics of flags associated to $C(h(x), \mathcal{H}_x)$ for $x\in \widetilde{S_g}$.
\end{proof}

We have also the following control of the position of the limit map compared to the conics in the fibration:

\begin{prop}
Let $\rho$ be a representation associated with a Higgs bundle in $\mathcal{S}_X^\R$. Let $C\subset \mathbb{RP}^2$ be any of the conics that define the fibration  of $\Omega_\rho$. The Anosov limit map $(\xi^1,\xi^2):\partial \Gamma_g\to \mathcal{F}_{\Delta}$ associated with $\rho$ is such that $\xi^1(x)$ is outside of $C$ and $\xi^2(x)$ intersects the interior of $C$ for all $x\in \partial \Gamma$.
\end{prop}

\begin{proof}
We know that $C$ lies in the domain of discontinuity $\Omega_\rho\subset \mathcal{F}_{1,2}$. In particular for all $y\in \partial \Gamma$, $\xi^1(y)$ does not lie on the conic $C$. Suppose that $\xi^1(\partial \Gamma)$ lies entirely in the inside of the conic $C$, then the corresponding loop in $\mathbb{RP}^2$ is contractible. However any Anosov representation of a closed surface group  with contractible first boundary map in $\SL(3,\R)$ must preserve a properly convex domain, and hence is necessarily Hitchin. 

The representations we consider are not in the Hitchin component, hence $\xi^1(\partial \Gamma)$ lies entirely in the outside of the conic $C$. The same argument in the dual of $\R^3$ implies the analog result for $\xi^2$.
\end{proof}


\newpage
\bibliographystyle{alpha}
\bibliography{references}
\end{document}